\documentclass[12pt,a4paper]{article}
\usepackage{amsfonts}
\usepackage[all]{xy}
\usepackage{graphics}

\newtheorem{theo}{Theorem}[section]

\newtheorem{coro}[theo]{Corollary}

\newtheorem{lm}[theo]{Lemma}

\newtheorem{prop}[theo]{Proposition}
\newtheorem{rem}[theo]{Remark}

\newtheorem{defi}{Definice}[section]
\newtheorem{de}{Definition}[section]
\newtheorem{kon}[defi]{Konstrukce}

\newtheorem{examp}{Example}[section]

\newenvironment{proo}{{\bf {Proof:}}}{\hfill $\square$}

\def\A{{\cal A}}

\def\C{{\cal C}}
\def\D{{\cal D}}
\def\E{{\cal E}}

\def\I{{\cal I}}

\def\NN{{\cal N}}
\def\M{{\cal M}}
\def\P{{\cal P}}

\def\S{{\cal S}}

\def\V{{\cal V}}

\def\Z{{\cal Z}}
\def\N{{\mathbb N}}

\def\KK{{\bf K}}
\def\MM{{\bf M}}

\def\alg{{\bf alg}}
\def\aalg{\!-\!{\bf alg}}
\def\Alg{{\bf Alg}\:}
\def\End{{\bf End}}

\def\ECl{\End_{\lambda}\C}

\def\colim{\mathop {\mbox{\rm colim}}\limits}
\def\copr{\displaystyle\coprod}
\def\pr{\displaystyle\prod}
\def\Id{{\rm Id}}

\def\id{{\rm id}}
\def\max{{\rm max}}
\def\Monad{{\bf Monad}}
\def\set{{\cal S}et}
\def\Ord{{\cal O}rd}

\def\ob{{\rm Ob}}
\def\mor{{\rm Mor}}
\def\iso{{\rm Iso}}
\def\imp{\Rightarrow}
\def\iff{\Leftrightarrow}
\def\hom{{\rm hom}}

\begin{document}
\begin{center}
\LARGE{\bf Varieties defined by natural transformations}\\
\vspace{1 cm}
\normalsize{Jan Pavl\'{\i}k} \\
\vspace{0.5 cm}
\small Faculty of Mechanical Engineering, Brno University of Technology,\\
Brno, Czech Republic\\
 {\tt pavlik@fme.vutbr.cz}\\
\end{center}

\begin{abstract} We define varieties of algebras for an arbitrary en\-do\-functor on a cocomplete category
using pairs of natural transformations. This approach is proved to be equivalent to the one of
equational classes defined by equation arrows. Free algebras in the varieties are investigated and
their existence is proved under the assumptions of accessibility.
\end{abstract}

In universal algebra we deal with varieties - classes of algebras satisfying a certain collection
of identities (pairs of terms of corresponding language). This concept was generalized by
Ad\'{a}mek and Porst in \cite{AP}. They worked with algebras for an endofunctor on
a cocomplete category and used the free-algebra construction developed by Ad\'{a}mek in \cite{A74}
(a chain of term-functors) to define the equation arrows as certain regular epimorphisms. Using
them they defined equational categories as analogy to varieties. These categories are later studied
in \cite{AT}.

We focus on another approach to varieties of algebras for a functor. We also use the free-algebra
construction and we define natural term as a natural transformation with codomain in a
term-functor. A pair of natural terms with the common domain will be called natural identity and
will be satisfied on an algebra if its both natural transformations have the same term-evaluation
on this algebra. Natural identities induce classes of algebras, which are proved to be precisely
the classes defined by means of equation arrows. We present several examples of such classes and
show that in some cases this approach essentially simplifies the presentation.

In the second chapter we investigate free algebras in a variety. The induction by natural
identities allows us to make the restriction on the identities with the domain preserving the
colimits of some small chains. Such identities will be called accessible. These cases still cover most of the usual examples and we prove that such varieties have free algebras. The proof uses the conversion of variety to a category of algebras for a diagram of monads used by Kelly in \cite{Kelly} to define algebraic colimit of
monads. His theorem proving the existence of free objects of this category yields the existence of
free algebras in the variety induced by accessible identities.
\paragraph{Notational convention}
The constant functor mapping the objects to object $X$ will be denoted by $C_X$. The initial object in a cocomplete category will be denoted by $0$. For functors, we omit the brackets and the composition mark $\circ$, when possible.
The class of objects and morphisms of a category will be denoted by $\ob$ and $\mor$, respectively.
The concrete isomorphism between two concrete categories over $\C$ (i.e. the isomorphism preserving the forgetful functor) will be denoted by $\cong_{\C}$.

\section{Classes of algebras}
\subsection{Algebras and equational classes}
\begin{de}
Let $F$ be an endofunctor on a category $\C$. By $\Alg F$ we denote the category of $F$-algebras -
its objects are pairs $(A, \alpha)$, where $\alpha:FA \to A$ is a morphism in $\C$. The morphism
$\phi_F:(A, \alpha) \to(B, \beta)$ of $F$-algebras is a morphism $\phi:A \to B$ such that $\phi
\circ \alpha = \beta \circ F\phi$. The subscript $_F$ in the notation of morphism is usually
omitted.\end{de}
\begin{rem}
There is a forgetful functor $\Z_F:\Alg F \to \C$ assigning to an algebra $(A,\alpha)$ its
underlying object $A$. \end{rem}
From now on, let $\C$ be a cocomplete category.
Let us recall the {\em free-algebra construction} (introduced in
\cite{A74}, generalized in \cite{AP}). We show the definition in the functorial form.
\begin{de}
By transfinite induction we define {\em term functors} $F_n:\C \to \C$, for $n \in \Ord$ and natural
transformations $w_{m,n}:F_m \to F_n$, for $m \leq n$:
\begin{itemize}
\item[] {\rm \underline{Initial step}:} $F_0=\Id_{\C}$, $w_{0,0}= \id$
\item[] {\rm \underline{Isolated step}:} Let $F_{n + 1}=FF_n + \Id_{\C}$,
the transformations $w_{0,n+1}= \iota_{n+1}$ and $q_n: FF_n \to F_{n+1}$ are the canonical injections
of $\Id_{\C}$ and  $FF_n$, respectively, into the coproduct and $w_{m+1,n+1}=[Fw_{m,n},\id_{\Id_{\C}}]$ for $m \leq n$ is defined by
$$\xymatrix{
FF_{m} \ar[d]_{q_{m}} \ar[rr]^{Fw_{m,n}}& &FF_n \ar[d]^{q_{n}}\\
F_{m+1} \ar[rr]^(.4){w_{m+1,n+1}}& &FF_{n}+\Id_{\C} = F_{n+1}\\
\Id_{\C} \ar[u]^{\iota_{m+1}}  \ar[urr]_{\iota_{n+1}}& & }.$$ If $m$ is a limit ordinal, then we
define $w_{m,m+1}$ as the unique factorization of $\{w_{k,m+1}| k < m\}$ over the colimit cocone
$\{w_{k,m}| k < m\}$.
\item[] {\rm \underline{Limit step}:} $F_n=\colim_{m < n} F_m$ and $w_{m,n}$ is the corresponding component of
colimit cocone.
\end{itemize}
The construction gives rise to the transformation $y_{n}:F \to F_n$ defined by $$y_{n}=w_{1,n} \circ q_0$$
for every ordinal $n>0$. 
\end{de}
To distinguish the transformations for different functors we denote the name of the functor in the
superscript: $w^{F}_{m,n}$, $q^{F}_{n}$, $\iota^{F}_{n}$, $y^{F}_{n}$.

The construction yields for every $m\leq n$ the property:
\begin{eqnarray}
\label{qw} w_{m,n}\circ  q_{m}  &=& q_{n} \circ Fw_{m,n}.
\end{eqnarray}
\begin{rem}\label{wv}
If we substitute $\Id_{C}$ by $C_0$ in initial step of construction, then we get equivalent concept.
\end{rem}
As a consequence of the definition we get the following properties (see \cite{AP}).
\begin{rem}\label{prop}
Given an $F$-algebra $(A, \alpha)$, for every $n \in \Ord$ there is a morphism (a {\em
term-evaluation} on $(A, \alpha)$) $$\epsilon_{n,(A, \alpha)}:F_nA \to A$$ defined recursively by:
$\epsilon_{0,(A,\alpha)} =\id_A$, $\epsilon_{n+1,(A,\alpha)}=[\alpha \circ
F\epsilon_{n,(A,\alpha)}, \id_A]$,
$$\xymatrix{
FF_n(A) \ar[d]_{q_{n,A}} \ar[rr]^{F\epsilon_{n,(A,\alpha)}}& &F(A) \ar[d]^{\alpha}\\
F_{n+1}(A) \ar[rr]^{\epsilon_{n+1,(A,\alpha)}}& &A\\
A \ar[u]^{\iota_{n+1,A}}  \ar[urr]_{\id_A}& & }$$ and by
$\epsilon_{l,(A,\alpha)}=\colim_{m<l}\epsilon_{m,(A,\alpha)}$ for a limit ordinal $l$. Then for
every $n$, $m \leq n$ we have:
\begin{eqnarray}
\label{epq}\epsilon_{n+1,(A,\alpha)} \circ q_{n,A}&=& \alpha \circ F\epsilon_{n,(A,\alpha)}\\
\label{epi}\epsilon_{n, (A, \alpha)} \circ \iota_{n, A}&=& \id_A\\
\label{epw} \epsilon_{m,(A, \alpha)} &=& \epsilon_{n,(A,\alpha)} \circ w_{m,n,A}\\
\label{epy} \epsilon_{n,(A, \alpha)} \circ y_{n,A}&=& \alpha,
\end{eqnarray}
where the last property requires $n >0$. We write the name of the functor in the superscript $\epsilon_{k,(A,\alpha)}=\epsilon^F_{k,(A,\alpha)}$, if necessary.
\end{rem}
We recall here the notion of equational category of $F$-algebras introduced in \cite{AP}.
\begin{de}
Let $X$ be an object of $\C$, $n \in \Ord$ . An {\em equation arrow of arity $n$ over $X$} is defined
as a regular epimorphism $e:F_nX \to E$. The object $X$ is called a {\em variable-object} of $e$.

We say, that an $F$-algebra $(A, \alpha)$ {\em satisfies} an equation arrow $e:F_nX \to E$ if for
every $f:X \to A$ there is a morphism $h:E \to A$ such that $\epsilon_{n,(A, \alpha)} \circ F_n f= h \circ e$.

For a class $\E$ of equation arrows we define an {\em equational class of $F$-algebras induced by
$\E$} as the class of all algebras satisfying all equations $e \in \E$. Considered as a full
subcategory of $\Alg F$ it is called an {\em equational category} and denoted by $\Alg (F,\E)$. If
$\E$ is a singleton, we say $\Alg (F,\E)$ is single-based.
%
\end{de}
As shown in \cite{AP}, this approach generalizes classical universal algebra on sets, since every
identity uniquely determines the regular epimorphism on the set of all terms, which is given by
unification of terms included in identity.

\subsection{Naturally Induced Classes}
Now we introduce the concept of algebras induced by natural transformations.
\begin{de}
Let $n$ be an ordinal and $G$ be a $\C$-endofunctor. A natural transformation $\phi:G \to F_n$ is
called a {\em natural term}, more precisely a {\em $G$-term}.

By {\em $G$-identity} we mean a pair of $G$-terms. Such pairs are called {\em natural identities}.
Let $\phi$ and $\psi$ be $m$-ary and $n$-ary $G$-terms, respectively. The functor $G$ is called a
{\em domain} and $(m,n)$ is an {\em arity-couple} of identity $(\phi,\psi)$. If $m=n$, we say
$(\phi,\psi)$ has the {\em arity} $n$.

We say, that an $F$-algebra $(A, \alpha)$ {\em satisfies} the $G$-identity $(\phi,\psi)$, if
$$\epsilon_{m, (A,\alpha)} \circ \phi_A = \epsilon_{n, (A,\alpha)} \circ \psi_A.$$
Then we write $$(A, \alpha)\models (\phi,\psi).$$

For a class $\I$ of natural identities we define a {\em naturally induced class of $F$-algebras} as
the class of all algebras satisfying all identities $(\phi,\psi)\in \I$. The corresponding full
subcategory of $\Alg F$ is denoted by $\Alg (F,\I)$. If $\I$ is a singleton, we say $\Alg (F,\I)$ is
single-induced.

Two natural identities are said to be algebraically equivalent iff they induce the same classes of $F$-algebras. Analogously we define the algebraic equivalence of classes of natural identities. The algebraic equivalence relation will be denoted by $\approx$.
\end{de}
\begin{rem}\label{simply}
\begin{enumerate}
\item Arities of components of a natural identity can be arbitrarily raised. Clearly for an identity $(\phi,\psi)$ of arity-couple $(m_1,m_2)$ we have $(\phi,\psi) \approx (w_{m_1,n} \circ \phi,w_{m_2,n} \circ \psi)$ for every $n \geq \max\{m_1,m_2\}$. Hence every natural identity is algebraically equivalent to the one consisting of natural terms of the same arity.
\item Every set $\NN=\{(\phi_i, \psi_i)|i \in I\}$ of $n$-ary natural identities is algebraically
 equivalent to the singleton. Clearly $\NN \approx \{(\phi,\psi)\}$, where $\phi, \psi$ are the unique factorizations of the cocones $\phi_i$, $\psi_i$, respectively, over the coproduct of domains of single identities.
\item As a consequence, every class naturally induced by a set of identities is single-induced.
\end{enumerate}
\end{rem}
\subsection{Examples of naturally induced classes}
In section \ref{conversion} we show, that every equational class is naturally
induced and vice versa. And since the concept of equational classes generalizes the varieties in universal
algebra, every variety of algebras in the classical sense is naturally induced class. Explicit correspondence is shown in the following example.
\begin{examp}
Let $\C=\set$, $\Sigma$ be a signature consisting of operation symbols $\sigma$ of (possibly infinite) arities
$ar(\sigma)$. Let $F=\coprod_{\sigma \in \Sigma} \hom(ar(\sigma),-)$ and $u_{\sigma}:\hom(ar(\sigma),-) \to F$ be the canonical inclusion for every $\sigma \in \Sigma$. Then the category of $\Sigma$-algebras is isomorphic to $\Alg F$. For each $\Sigma$-term $\tau$ let $X_{\tau}$ be the set of variables occurring in $\tau$ and let $d(\tau)$ be the depth of $\tau$ (supremum of ordinals corresponding to chains of the proper subterms of $\tau$ ordered by subterm-relation).

For a given term $\tau$ we assign a $d(\tau)$-ary $G_{\tau}$-term $\phi_{\tau}$, where $G_{\tau}=\hom(X_{\tau},-)$. The transformation $\phi_{\tau}$ is defined inductively: if $\tau$ is a variable $x$, then $\phi_{\tau}:\hom(\{x\},-) \to F_0$ is obvious isomorphism. If $\tau=\sigma(\rho_{i};i \in
ar(\sigma))$ and we have $\phi_{\rho_{i}}:\hom(ar(\rho_{i}),-) \to F_{d(\rho_{i})}$ for each $i
\in ar(\sigma)$, then we can extend all transformations $\phi_{\rho_{i}}$ to
$\phi'_i:\hom(ar(\rho_{i}),-) \to F_{n}$, where $n=\sup\{d(\rho_i)|i\in ar(\sigma)\}$. We define
$\phi_{\tau}$ the following way. Since $X_{\rho_i} \subseteq X_{\tau}$ for every $i$, we have
$p_i:\hom(X_{\tau},-) \to \hom(X_{\rho_i},-)$, hence the factorization over the limit cone yields
the unique $r:\hom(X_{\tau},-) \to \prod_{i \in ar(\sigma)} \hom(X_{\rho_i},-)$. We define $\phi_{\tau}$ as the following composition:
$$\xymatrix{
\hom(X_{\tau},-) \ar[d]^{\phi_{\tau}} \ar[r]^(.4){r} & \pr_{i \in ar(\sigma)}\hom(X_{\rho_i},-) \ar[r]^(.6){\prod \phi'_i}
&\pr_{i \in ar(\sigma)}F_n \ar@{=}[d]^{iso}\\
F_{n+1}&FF_n \ar[l]_{q_{n}}&\hom(ar(\sigma),-) \circ F_n \ar[l]_(0.6){u_{\sigma}F_n}
}$$ Observe, that $n+1=d(\tau)$.

For each $\Sigma$-term we have assigned a natural term. Now for identity $(\tau_1,\tau_2)$ consisting of two $\Sigma$-terms with variables in $X$ we assign a pair of corresponding natural $\hom(X,-)$-terms. It is easy to see, that we get the identity, which induces exactly the variety given by $(\tau_1,\tau_2)$. For example, monoids are objects of $\Alg((\hom(2,-)+\hom(0,-), \{i,j,k\})$, where $i$ is binary identity with domain $\hom(3,-)$ and stands for associativity while $j,k$ are unary with domain $\Id$ and correspond to left and right neutrality of $1$.
\end{examp}
The concept can be used to define naturally induced classes of algebras even on some illegitimate
categories.
\begin{examp}\label{monads}
Let $\C= \End \A$ be the illegitimate category of endofunctors on some cocomplete category $\A$.
The composability of objects of $\C$ yields for every $k \in \omega$ the existence of the "composition power functor" $S_k:\C \to \C$  such that $$S_k(P)=\underbrace{P\circ P \circ \dots \circ P}_{k\: \rm times}.$$
We can define analogies for universal algebras - all we need to do, is to substitute products of sets by composition of functors and each $\hom(k,-)$ by $S_k$ in the description above. As analogy of monoids we get the category $\Monad \:\A$ of monads on $\A$. Namely, $\Monad \:\A = \Alg ((S_2+S_0),
\{i,j,k\})$, where domains of identities $i,j,k$ are $S_3,S_1,S_1$, respectively. Each operation
$\pi:(S_2+S_0)(P)\to P$ decomposes into $\mu:S_2(P)=PP\to P$ and $\eta:S_0=\Id \to P$ and the
satisfaction of identities fully corresponds to usual condition claimed on $\mu$ and $\eta$.
\end{examp}
Theorem 3.6 in \cite{AP} describes the equational presentation of category of algebras for a monad. The following example shows their the presentation by natural identities.
\begin{examp}\label{monad}
Given a monad ${\bf M}= (M, \eta, \mu)$ on $\C$, then its Eilenberg-Moore category $\MM\aalg$ is
a class of $M$-algebras induced by two natural identities:
$$\xymatrix{
M^2 \ar[dr]_{Mq_0} \ar[r]^{\mu} &M \ar[r]^{q_0}&M_1\\
& MM_1 \ar[r]^{q_1}&M_2,}\quad \xymatrix{
\Id \ar[d]^{\eta} \ar[r]^{\id}& M_0\\
M \ar[r]^{q_0}&M_1.} $$Therefore
$$\MM\aalg=\Alg (M, \{(q_0 \circ \eta, \id_{\Id}), (q_1 \circ Mq_0, q_0 \circ \mu)\}.$$
\end{examp}
\begin{examp}\label{powerset}
Consider the power-set monad on $\set$ defined by power-set functor $\P$ and transformation $\eta:
\Id_{\set} \to \P$, $\mu:\P^2 \to \P$ given by assignments $\eta_X(x)=\{x\}$, $\mu_X(\{X_i|i\in
I\})=\bigcup_{i \in I}X_i$. As a concrete instance of the previous case for power-set monad $(\P,\eta,
\mu)$ we get the category of join-complete semilattices {\bf JCSlat}. Hence this class is
presentable by a pair of naturally induced identities - compare with presentation by a proper class
of equation arrows (see \cite{AP}, Example 3.3 - we need equation arrows $e_X:F_3X \to E_X$ for
every set $X$).
\end{examp}
\subsection{Conversion theorem}\label{conversion}
Our aim is to prove that naturally induced classes and equational classes coincide. At first we show that every single-based
equational class is naturally induced. Then, conversely, we prove that every class induced by a single natural identity is equational. The crucial point of the proof is the local smallness of category $\C$.
\begin{rem}
Within the proof we use the copower functor:\\
Given an object $Q \in \C$, there is a functor $- \bullet Q: \set \to
\C$, which is left adjoint to $\hom(Q,-):\C \to \set$. It assigns to a set $M$ the coproduct of $M$ copies of $Q$ (the
"$M$-th" copower of $Q$) and for a mapping $h:M \to N$ we define $h \bullet Q$ as the unique factorization of cocone
$u_{h(m)}:Q \to \copr_{j \in N}Q, m \in M,$ over a colimit cocone $u_{m}:Q \to \copr_{j \in M}Q$.

Then we get the adjunction $(\eta,\varepsilon):(- \bullet Q) \dashv \hom(Q,-):\C \to \set$, where the unit morphism $\eta_{X}:X \to
\hom(Q,X \bullet Q)$ for a set $X$ and $x \in X$ is defined by $\eta_{X}(x)=u_{x}:Q \to X \bullet Q$, i.e. the
$x$-labeled canonical injection into the coproduct. Moreover, for an object $A$ of $\C$, the counit $\varepsilon:
\hom(Q,A) \bullet Q \to A$ is defined as the unique factorization of a cocone $\{\phi|\phi:Q \to A\}$ over the colimit.
\end{rem}
\begin{lm}\label{singvar}
Every single-based equational class is naturally single-induced class.
\end{lm}
\begin{proo}
Let $\S$ be a single-based equational class of $F$-algebras defined by an equation arrow $e$, where $e$ is a regular epimorphism $F_n X \to E$ such that $(E,e)$ is a coequalizer of $\phi_0,\psi_0:\xymatrix{
Q \ar@<1ex>[r]\ar@<-1ex>[r]&F_n X}$.
We define a mapping $\theta_{\phi,A}:\hom(X,A) \to \hom(Q,F_n A)$. For every
$f:X \to A$ let $\theta_{\phi,A}(f) = F_n f \circ \chi_0:Q \to F_n A$. Now let
$$G= (- \bullet Q) \circ \hom(X,-),$$
$$\phi_A = \widetilde{\theta_{\phi,A}}: \hom(X,A)\bullet Q \to F_n A.$$
Clearly $\phi_A$ is the component of a natural transformation $\phi:G\to F_n$.
Observe, that for every $f:X \to A$, holds $$\phi_A \circ u_f =\theta_{\phi,A}(f)= F_nf \circ
\phi_0.$$
Analogously we define the natural transformations $\theta_{\psi,-}:G \to F_n$ and
$\psi:G \to F_n$ satisfying $\psi_A \circ u_f =F_nf \circ
\psi_0$.
Now we have the functor $G$ and $G$-identity $(\phi, \psi)$. It remains to show, that it induces exactly the equational class $\S$.

Let $(A, \alpha)$ satisfy the equation arrow $e$. Then for every $f:X \to A$ there is $h:E \to A$,
such that $\epsilon_{n, (A, \alpha)} \circ F_n f = h \circ e$. Then we have
$$\begin{array}{rcl}
\epsilon_{n, (A, \alpha)} \circ \phi_A \circ u_f&=&\epsilon_{n, (A, \alpha)} \circ F_nf \circ \phi_0\\
&=&h \circ e \circ \phi_0 = h \circ e \circ \psi_0 \end{array}$$ and by symmetry we get $\epsilon_{n, (A, \alpha)}
\circ \phi_A \circ u_f=\epsilon_{n, (A, \alpha)} \circ \psi_A \circ u_f$. Since $f$ was chosen arbitrarily and injections $u_f$
form a colimit cocone, we have $\epsilon_{n, (A, \alpha)} \circ \phi_A=\epsilon_{n, (A, \alpha)} \circ \psi_A$, i.e.
$(A, \alpha)$ satisfies the $G$-identity $(\phi,\psi)$.

Now let $(B,\beta)$ be an $F$-algebra in a class induced by the $G$-identity $(\phi,\psi)$. Let $g:X \to B$ be a morphism in $\C$. Then we
have
$$\begin{array}{rcl} \epsilon_{n, (B,\beta)} \circ F_n g \circ \phi_0&=& \epsilon_{n, (B,\beta)} \circ \phi_B
\circ u_g\\
&=&\epsilon_{n, (B,\beta)} \circ \psi_B \circ u_g
\end{array}$$
and again by symmetry we get $\epsilon_{n, (B,\beta)} \circ F_n g \circ \phi_0 = \epsilon_{n,
(B,\beta)} \circ F_n g \circ \psi_0$, hence $\epsilon_{n, (B,\beta)} \circ F_n g$ coequalizes the
pair $(\phi_0, \psi_0)$ and there is a unique $h:E \to B$ such that $\epsilon_{n, (B,\beta)} \circ
F_n g = h\circ e$. Thus $(B, \beta)$ satisfies the equation arrow $e$.
\end{proo}
\begin{lm}\label{singnat}
Every naturally single-induced class is equational.
\end{lm}
\begin{proo}
Let $G$ be a $\C$-endofunctor. Let $\NN$ be a class induced by a $G$-identity $(\phi,\psi)$. Due
to Remark \ref{simply} we may assume that $\phi$ and $\psi$ have the same arity, say $n$, therefore
both are the natural transformations $G \to F_n$. Let $(E,e)$ be the coequalizer of $\phi$ and
$\psi$. Then for every object $X$ of $\C$ we have a morphism $e_X:F_nX \to EX$. Let $\E=\{e_X|X \in
\ob\C\}$. We will prove $\NN =\Alg (F,\E)$.

Let $(A,\alpha)$ satisfy $(\phi,\psi)$. Then for every $X \in \ob\C$ and $f:X \to A$ we have
$$\begin{array}{rcl}
\epsilon_{n, (A,\alpha)} \circ F_nf \circ \phi_X &=&\epsilon_{n, (A,\alpha)} \circ \phi_A \circ Gf\\
 &=&\epsilon_{n, (A,\alpha)} \circ \psi_A \circ Gf\\
 &=&\epsilon_{n, (A,\alpha)} \circ F_nf \circ \psi_X,
 \end{array}$$ therefore we have coequalizing morphism $\epsilon_{n, (A,\alpha)} \circ F_nf$ for $(\phi_X, \psi_X)$.
Since the colimits of functors are calculated componentwise, $e_X$ is a coequalizer of $(\phi_X,
\psi_X)$, therefore there is a unique $h:EX \to A$, such that $\epsilon_{n, (A,\alpha)} \circ F_nf
= h \circ e_X$.

Given an $F$-algebra $(B, \beta)$ satisfying all equation arrows from $\E$, then it satisfies the
arrow $e_B:F_nB \to EB$ and there is $h:EB \to B$ (chosen for $\id_B:B \to B$) such that
$\epsilon_{n,(B, \beta)} = h \circ e_B$, therefore the property is satisfied, since $e_B$
coequalizes the pair $(\phi_B, \psi_B)$.
\end{proo}
\begin{theo}\label{nat}
Let $F$ be an endofunctor on a cocomplete category $\C$. Then the equational classes of $F$-algebras coincide with naturally induced classes of $F$-algebras.
\end{theo}
\begin{proo}
Every equational class $\S$ is a (possibly large) intersection of single-based ones and those are by the Lemma \ref{singvar} naturally induced, hence we get $\S$ to be naturally induced by the class of all natural identities present in some of the collection. Conversely, the naturally induced class $\NN$ is a (possibly large) intersection of the ones induced by a single natural identity, which are due to Lemma \ref{singnat} equational classes induced by a class of equation arrows. Union of these classes defines the class of all equation arrows defining the class $\NN$ as an equational class.
\end{proo}
\begin{de}
Class of algebras induced by equations or natural identities is called a {\em variety}.
\end{de}
\section{Free algebras in the variety}\label{free}
Our aim is to answer the question of existence of free algebras in varieties. At first we recall well-known results involving free algebras, which will be used to solve this question.
\subsection{Free algebras and monads}
We will work with $\C$-endofunctors preserving the colimits of $\lambda$-labeled chains, where $\lambda$ is infinite limit ordinal - let the class containing these functors be denoted by $\ECl$. Since colimits commute with colimits we get the following property (see also \cite{Kelly}, 2.4.).
\begin{prop}\label{close}
The class $\ECl$ is closed under colimits and compositions.
\end{prop}
\begin{de}
Let $G$ be an endofunctor on $\C$. The natural $G$-identity is called {\em accessible} if $G$ preserves the colimits of $\lambda$-labeled chains for some infinite limit ordinal $\lambda$.
\end{de}
\begin{de}
The functor $F$ admitting free $F$-algebras is called a {\em varietor}.
\end{de}

Let $F$ preserve the colimits of $\lambda$-indexed chains. Then, as shown in \cite{A74}, $F$ is a varietor. Since $F$
preserves the $\lambda$-labeled chains, $FF_{\lambda}$ is a colimit of chain
 $\{FF_{n}|n < \lambda\}$. Hence one can see that $w_{\lambda,\lambda +1}$ is isomorphism.
 In such a case we say that {\em the free $F$-algebra construction stops after $\lambda$ steps}. If we set
 $\upsilon = \colim_{n<\lambda} q_n$, on every object $A$ we get the free $F$-algebra $$\V_F(A)=(F_{\lambda}A, \upsilon_A).$$
 If necessary, we write the name of functor $F$ in the superscript: $\upsilon=\upsilon^F$.

This construction gives rise to the functor $\V_F:\C \to \Alg F$, $\V_F=(F_{\lambda},\upsilon)$ together with
transformation $\epsilon_{\lambda}:\V_F\Z_F \to \Id_{\Alg F}$, $\epsilon_{\lambda,(A,\alpha)}:(F_{\lambda}A,\upsilon_A)
 \to (A,\alpha)$. Hence we have got the free functor $\V_F$ and adjunction $\V_F \dashv \Z_F$, where the unit and counit
 are $\iota_{\lambda}$ and $\epsilon_{\lambda}$, respectively.

It is well known fact, that this adjunction yields the free monad over a functor $F$ - see
\cite{AHS}, Theorem 20.56. Hence the free monad over $F$ is $$\M(F)=(F_{\lambda},\eta^F, \mu^F),$$ where
$\eta^F=\iota^F_{\lambda}$ and $\mu^F=\Z_F\epsilon_{\lambda}\V_F$ and universal morphism for $F$ is $y^F_{\lambda}:F \to F_{\lambda}$. More detailed approach to the theory
of monads can be found in \cite{AHS}, \cite{Barr} and \cite{Kelly}.

Now we use another functor $G$ in $\ECl$ and we work with its algebras. We point out the important consequences of the discussion above:
\begin{prop}\label{monadj}
Let there be a transformation $\rho:G \to F_{\lambda}$. Then there is a transformation $\sigma:G_{\lambda}\to
F_{\lambda}$, subject to the conditions:
\begin{enumerate}
\item $\sigma=\overline{\rho}$ is given by the freeness of $\M(F)$ as the unique monad transformation $\M(G) \to \M(F)$ corresponding to $\rho:G \to F_{\lambda}$; thus
    $$\sigma \circ y^G_{\lambda}=\rho.$$
\item $\sigma_A=\widetilde{\eta^F_A}$ is given by the adjunction $\V_G \dashv \Z_G$ as the unique $G$-algebra morphism $\V_G(A) \to P(A)$ corresponding to $\eta^F_A:A \to F_{\lambda}A=\Z_G P(A)$, where $P:\C \to \Alg G$ is the functor assigning to an object $A$ an algebra $(F_{\lambda}A,\beta_A)$ and $\beta_A = (\mu^F \circ \rho F_{\lambda})_A$. Hence
    $$\sigma \circ \iota^G_{\lambda}=\eta^F.$$
\item For $k\leq \lambda $ let $\epsilon^G_{k,P}:G_kF_{\lambda} \to F_{\lambda}$ be the obvious transformation with the components $\epsilon^G_{\lambda, P(A)}$. Then the following property holds:
    $$\sigma = \epsilon^G_{\lambda,P} \circ G_{\lambda}\eta^F.$$
\end{enumerate}
\end{prop}

We will show that $\sigma$ defined above from the transformation $\rho:G \to F_{\lambda}$ can be gained via the colimit construction, which will be later useful.

\begin{de}\label{rokon}
For all $k \in \N$ we define the transformations $\rho_k:G_k \to F_{\lambda}$ inductively:
$\rho_1=[\rho,\eta^F]$, $\rho_{k+1}=[\mu^F \circ \rho F_{\lambda} \circ G\rho_k,\eta^F]$
\end{de}
\begin{lm}
For every $j<k\in \N$ holds $\rho_k \circ w^{G}_{j,k}=\rho_j$.
\end{lm}
\begin{proo}
For every natural $j<k$ we have
$$\begin{array}{rcl}
\rho_k \circ w^{G}_{j,k} &=&[\mu^F \circ \rho F_{\lambda} \circ G\rho_k,\eta^F]\circ w^{G}_{j,k}\\
&=& [\mu^F \circ \rho F_{\lambda} \circ G\rho_{k-1} \circ Gw^{G}_{j-1,k-1},\eta^F]\\
&=&[\mu^F \circ \rho F_{\lambda} \circ G(\rho_{k-1} \circ w^{G}_{j-1,k-1}),\eta^F]
\end{array}$$
If $j=1$ then holds
$$\begin{array}{rcl}
\mu^F \circ \rho F_{\lambda} \circ G(\rho_{k-1} \circ w^{G}_{j-1,k-1})&=&
\mu^F \circ \rho F_{\lambda} \circ G(\rho_{k-1} \circ w^{G}_{0,k-1})\\
&=& \mu^F \circ \rho F_{\lambda} \circ G(\rho_{k-1} \circ \iota^{G}_{k-1})\\
&=& \mu^F \circ \rho F_{\lambda} \circ G\eta^{F}\\
&=& \mu^F \circ F_{\lambda}\eta^{F} \circ \rho = \rho
\end{array}$$
hence the property holds for $j=1$ and every $k>1$.

Now let $1<j<k$ and assume the validity of $\rho_{k-1} \circ w^{G}_{j-1,k-1}=\rho_{j-1}$. Then we have: $\mu^F \circ \rho F_{\lambda} \circ G(\rho_{k-1} \circ w^{G}_{j-1,k-1})=\mu^F \circ \rho
F_{\lambda} \circ G\rho_{j-1}$ hence $\rho_k \circ w^{G}_{j,k}=\rho_j$ and the proof is complete.
\end{proo}

Now the transformations $\rho_k$ form a compatible cocone for $w^{G}_{j,k}$ hence we can extend it
to the infinite limit case:
\begin{de}\label{rolim}
For a limit ordinal $l$ let $\rho_{l}=\colim_{k<l}\rho_{k}$.
\end{de}
Since the $w$-compatibility clearly holds, our construction extends to the ordinal chain with
analogous definition for the isolated step as for finite indices. The definition yields the property for every $k\leq \lambda$:
\begin{eqnarray}
\label{etag} \rho_{k}\circ \iota^G_{k}&=&\iota^F_{\lambda}
\end{eqnarray}

To prove that this chain of transformations converges to $\sigma$, we will show that its colimit
is $G$-algebra morphism.
\begin{lm}
The transformation $\rho_{\lambda}:G_{\lambda} \to F_{\lambda}$ underlies the natural
transformation $\V_G(A) \to P(A)$ of functor $\C \to \Alg G$, where $P$ is the functor used in Proposition \ref{monadj}.
\end{lm}
\begin{proo}
What we need to prove is that for an object $A$ in $\C$ the morphism $\rho_{\lambda,A}:G_{\lambda}A \to F_{\lambda}A$ is a
 $G$-algebra morphism $(G_{\lambda}A,\upsilon^G_A) \to (F_{\lambda}A,\beta_A)$. It suffices to prove the equality of natural transformations: $\beta \circ G\rho_{\lambda}=\rho_{\lambda} \circ \upsilon^G$.
Let $k<\lambda$, then we have
$$\begin{array}{rcl}
\rho_{\lambda} \circ \upsilon^G \circ Gw^{G}_{k, \lambda}&=&
\rho_{\lambda} \circ w^{G}_{k+1, \lambda}\circ q^G_{k} \\
&=& \rho_{k+1} \circ q^G_{k} \\
&=&\mu^F \circ \rho F_{\lambda}\circ G\rho_k\\
&=&\beta \circ G\rho_k\\
&=&\beta \circ G\rho_{\lambda} \circ Gw^{G}_{k, \lambda}
\end{array}$$
and since $G$ preserves the colimits of $\lambda$-indexed chains, $\{Gw^{G}_{k, l}|k\leq l<
\lambda\}$ is the colimit cocone. From the uniqueness of factorization over the colimit we get the
required equality.
\end{proo}
\begin{lm}\label{rosi}
The transformations $\rho_{\lambda}, \sigma:G_{\lambda} \to F_{\lambda}$ coincide.
\end{lm}
\begin{proo}
Given a $\C$-object $A$, then due to previous lemma, $\rho_{\lambda,A}$ is a
 $G$-algebra morphism $\rho_{\lambda,A}:(G_{\lambda}A,\upsilon^G_A) \to (F_{\lambda}A,\beta_A)$ and
 by (\ref{etag}) we have $\rho_{\lambda}\circ \iota^G_{\lambda}=\iota^F_{\lambda}=\eta^F$, hence by uniqueness of factorization
of $\eta^F_A:A \to \Z_G(F_{\lambda}A,\beta_A)$ over $\eta^G_{A}=\iota^G_{\lambda}$ we get
$\rho_{\lambda,A}=\widetilde{\eta^F_A}$ which is due to Proposition \ref{monadj} equal to
$\sigma_A$.
\end{proo}
\begin{lm}\label{equi}
Let $\phi,\psi:G \to F_{\lambda}$ be natural transformations. Then for every $k \leq \lambda$ we
have the algebraic equivalence:$$(\phi,\psi) \approx (\phi_{k},\psi_{k}),$$ where $\phi_k,\psi_k$
are derived from $\phi, \psi$, respectively, as in Definition \ref{rokon}, \ref{rolim}.
\end{lm}
\begin{proo}
Let  $(A,\alpha)$ be an $F$-algebra. Then for every $k \leq \lambda$
$$(\ast) \qquad\epsilon_{\lambda,(A,\alpha)} \circ \phi_{k,A} \circ \iota^G_k=\epsilon_{\lambda,(A,\alpha)} \circ \psi_{k,A} \circ \iota^G_k,$$
since $\epsilon_{\lambda,(A,\alpha)} \circ \phi_{k,A} \circ \iota^G_k=\epsilon_{\lambda,(A,\alpha)}
\circ \eta^F=\epsilon_{\lambda,(A,\alpha)} \circ w_{0,\lambda,A}=\epsilon_{0,(A,\alpha)} =\id$.

Let $(A,\alpha) \models (\phi,\psi)$, i.e.
$$(h) \qquad\epsilon_{\lambda,(A,\alpha)} \circ \phi_{A}=\epsilon_{\lambda,(A,\alpha)} \circ \psi_{A}.$$
We will show by induction, that then $(A,\alpha) \models (\phi_k,\psi_k)$ for every $k \leq \lambda$. In each step we shorten the computations using the $(\phi - \psi)$-symmetry of expressions.
\begin{itemize}
\item[]\underline{$k=1$}: \quad
Since $\epsilon_{\lambda,(A,\alpha)} \circ \phi_{1,A} \circ q^G_{0,A}=\epsilon_{\lambda,(A,\alpha)} \circ \phi_{A}$, then from $(h)$ and symmetry we get $\epsilon_{\lambda,(A,\alpha)} \circ \phi_{1,A} \circ q^G_{0,A}=\epsilon_{\lambda,(A,\alpha)} \circ \psi_{1,A} \circ q^G_{0,A}$ which together with $(\ast)$ yields $\epsilon_{\lambda,(A,\alpha)} \circ \phi_{1,A} =\epsilon_{\lambda,(A,\alpha)} \circ \psi_{1,A}$.
\item[]\underline{$1<k<\lambda$, $k$ isolated}: \quad Assume the hypothesis
$$(h_k) \qquad\epsilon_{\lambda,(A,\alpha)} \circ \phi_{k,A} =\epsilon_{\lambda,(A,\alpha)} \circ \psi_{k,A}$$
Recall, that $\epsilon_{\lambda,(A,\alpha)}:F_{\lambda}A \to A$ is a morphism $(F_{\lambda}A, \mu^F) \to (A, \epsilon_{\lambda,(A,\alpha)})$, i.e. holds $\epsilon_{\lambda,(A,\alpha)} \circ \mu^F = \epsilon_{\lambda,(A,\alpha)} \circ F_{\lambda}\epsilon_{\lambda,(A,\alpha)}$. Then we have:
$$\begin{array}{rcl}
\epsilon_{\lambda,(A,\alpha)} \circ \phi_{k+1,A} \circ q^G_{k,A}
&=&\epsilon_{\lambda,(A,\alpha)} \circ \mu^F_A \circ \phi_{F_{\lambda}A} \circ G\phi_{k,A}\\
&=&\epsilon_{\lambda,(A,\alpha)} \circ F_{\lambda}\epsilon_{\lambda,(A,\alpha)} \circ F_{\lambda}\phi_{k,A} \circ \phi_{G_kA}\\
&\stackrel{(h_k)}{=}&\epsilon_{\lambda,(A,\alpha)} \circ F_{\lambda}\epsilon_{\lambda,(A,\alpha)} \circ F_{\lambda}\psi_{k,A} \circ \phi_{G_kA}\\
&=&\epsilon_{\lambda,(A,\alpha)} \circ F_{\lambda}\epsilon_{\lambda,(A,\alpha)}
\circ \phi_{F_{\lambda}A} \circ G\psi_{k,A} \\
&=&\epsilon_{\lambda,(A,\alpha)}
\circ \phi_{A} \circ G\epsilon_{\lambda,(A,\alpha)}\circ G\psi_{k,A} \\
&\stackrel{(h)}{=}&\epsilon_{\lambda,(A,\alpha)}
\circ \psi_{A} \circ G\epsilon_{\lambda,(A,\alpha)}\circ G\psi_{k,A} \\
&\stackrel{symmetry}{=}&\epsilon_{\lambda,(A,\alpha)} \circ \psi_{k+1,A} \circ q^G_{k,A}
\end{array}$$
and together with $(\ast)$ we get $\epsilon_{\lambda,(A,\alpha)} \circ \phi_{k+1,A} =\epsilon_{\lambda,(A,\alpha)} \circ \psi_{k+1,A}$.
\item[]\underline{$l\leq\lambda$, $l$ limit}: \quad Assume $(h_k)$ for every $k<l$. Then from the uniqueness of factorization of cocone
$\epsilon_{\lambda,(A,\alpha)} \circ \phi_{k,A}:G_{k}\to A$ over the colimit cocone $w^G_{k,l,A}$
we get $\epsilon_{\lambda,(A,\alpha)} \circ \phi_{l,A}=\epsilon_{\lambda,(A,\alpha)} \circ
\psi_{l,A}$.
\end{itemize}
We have proved for every $k$: $$(A,\alpha)\models (\phi,\psi) \quad \imp \quad (A,\alpha)\models
(\phi_k,\psi_k).$$ However $(A,\alpha)\models (\phi_k,\psi_k)$ easily implies $(A,\alpha)\models
(\phi,\psi)$ since $\phi=\phi_{1} \circ q^G_{0}=\phi_{k}\circ w^G_{1,k} \circ q^G_{0}=\phi_{k}\circ
y^G_{k}$.
\end{proo}

\subsection{Algebras for a diagram of monads}
This section refers to the paper \cite{Kelly} of G. M. Kelly, chapter VIII., which deals with
colimits of monads. It is well known (see e.g. \cite{AHS}, Corollary 20.57), that for every
varietor $F$ the categories of its algebras and algebras for a monad $\M(F)$ are concretely
isomorphic via the comparison functor.

Let $D:\D \to \Monad \:\C$ be a diagram and $D(x)=(M_x,\eta^x,\mu^x)$ for every object $x \in \D$.
Consider the category $D\!-\!\alg$ of {\em algebras for a diagram $D$ of monads}, whose objects are
collections of $\C$-morphisms $\{\alpha_x:M_xA \to A|x \in \D\}$, where $\alpha_x$ is in
$D(x)\aalg$ for every object $x \in \D$ and for each $f:x \to y$ in $\D$ the $D$-compatibility
condition $\alpha_y \circ D(f)_A =\alpha_x$ is satisfied. The morphisms in $D\!-\!\alg$ are the
morphisms of algebras for each $x$, i.e. $\phi:(A,\alpha)\to (B, \beta)$ is a morphism if $\phi
\circ \alpha_x = \beta_x \circ M_x(\phi)$ for every $x$. If there is a monad $\KK$ such that $\KK\aalg \cong_{\C} D\!-\!\alg$, then this monad is called {\em algebraic colimit of
$D$}.

Kelly asked about existence of this algebraic colimit. It came out to be equivalent to existence
of the free objects in $D\!-\!\alg$. He proved the existence in his Theorem 27.1 in \cite{Kelly} under
the general assumptions of existence of suitable factorization systems and some smallness
requirements for the monads. Using the trivial factorization system $(\iso,\mor)$ and preservation of
$\lambda$-labeled chains, we get this theorem in the following form:
\begin{theo}\label{kelly}
Let underlying functor of each $D(x)$ preserve the colimits of $\lambda$-labeled chains. Then
$D\!-\!\alg$ has the free objects.
\end{theo}

This theorem will be used to prove the existence of free object in a variety induced by accessible identities. Let $F$ be a functor
in $\End_{\kappa}\C$ for some infinite limit ordinal $\kappa$ and consider the variety of
$F$-algebras induced by a set of accessible natural identities. Since the free $F$-algebra
construction stops after $\kappa$ steps, we may consider arity of each natural term to be less or
equal to $\kappa$. Then, due to the Remark \ref{simply}, the set of natural identities can be
substituted by a single identity $(\phi,\psi)$. Its domain, denoted by $G$, is the coproduct of
domains of single identities, hence, due to \ref{close}, it preserves colimits of $\nu$-indexed
chains for some large enough limit ordinal $\nu$. Let $\lambda = \max\{\kappa,\nu\}$, then $F,G \in \ECl$. Hence
the arity of $(\phi,\psi)$ can be set to $\lambda$.

Let $\D$ be a category consisting of two objects $0,1$, their identities and two more morphisms
$f,g:0 \to 1$. Let $D:\D \to \Monad \:\C$ be a diagram such that $D(0)=\M(G)$, $D(1)=\M(F)$,
$D(f)=\overline{\phi}$, $D(g)=\overline{\psi}$, where $\overline{\phi},\overline{\psi}$ are the monad transformations given by Proposition \ref{monadj}. We will prove the concrete equivalence of $\Alg(F,
(\phi,\psi))$ and $D\!-\!\alg$.

\begin{lm}
For the $\lambda$-ary $G$-identity $(\phi,\psi)$ and diagram $D$ defined above holds:$$\Alg(F,
(\phi,\psi)) \cong_{\C} D\!-\!\alg.$$
\end{lm}
\begin{proo}
Consider the comparison functor $I:\Alg F \to \M(F)\aalg$ assigning to an $F$-algebra $(A, \alpha)$
the $F_{\lambda}$-algebra $(A, \epsilon_{\lambda,(A,\alpha)})$. Then due to Lemma \ref{equi}
$$(A, \alpha) \models (\phi,\psi) \quad \iff \quad(A, \alpha) \models (\phi_{\lambda},\psi_{\lambda})\quad \iff
\quad \epsilon_{\lambda,(A,\alpha)} \circ \phi_{\lambda} = \epsilon_{\lambda,(A,\alpha)} \circ \phi_{\lambda}.$$
For every $F_{\lambda}$-algebra $(A,\beta)$ holds $\epsilon^{F_{\lambda}}_{1,(A,\beta)} \circ q^{F_{\lambda}}_0=\beta$, hence we get $$\epsilon_{\lambda,(A,\alpha)} \circ \phi_{\lambda} = \epsilon_{\lambda,(A,\alpha)} \circ \phi_{\lambda}\quad \iff \quad I(A, \alpha) \models (\phi^{\ast},\psi^{\ast}),$$ where $\phi^{\ast}=q^{F_{\lambda}}_0 \circ \phi_{\lambda}$.
Since $I$ is the isomorphism with an inverse given by $(A, \beta) \mapsto (A, \beta \circ
y_{\lambda,A})$, we get $$\Alg (F,(\phi,\psi)) \cong_{\C} \M(F)\aalg \cap \Alg
(F_{\lambda},(\phi^{\ast},\psi^{\ast})).$$ Due to Proposition \ref{monadj} we have
$\phi_{\lambda}=\overline{\phi}$, which is a monad transformation (and analogously for $\psi$), hence we
get $\M(F)\aalg \cap \Alg (F_{\lambda},(\phi^{\ast},\psi^{\ast}))$ to be concretely
isomorphic to $D\!-\!\alg$. 
\end{proo}

Now we can use Kelly's theorem to conclude our investigation:
\begin{theo}\label{freevar}
Let $F$ preserve the colimits of $\lambda$-indexed chains for some limit ordinal $\lambda$. Then the free algebra exists in every variety induced by a set of  accessible identities.
\end{theo}
To express the consequence for the varieties presented by equation arrows, recall the notion of presentability of an object (see \cite{AR}):
\begin{de}
Let $\lambda$ be a regular cardinal. An object $A$ of a category is called {\em $\lambda$-presentable} provided that its hom-functor $\hom(A,-)$ preserves $\lambda$-directed colimits. An object is called {\em presnetable} if it is $\lambda$-presentable for some $\lambda$.
\end{de}
\begin{coro}
Let $F$ preserve the colimits of $\lambda$-indexed chains for some limit ordinal $\lambda$. Then the free algebra exists in every variety induced by a set of equation arrows with presentable variable-objects.
\end{coro}
\begin{proo}
As shown in the proof of Lemma \ref{singvar}, an equation arrow $e:F_nX \to E$ converts to a natural identity with
the domain $G=(- \bullet Q) \circ \hom(X,-)$ for some $Q \in \ob\C$. If the variable-object $X$ is presentable, $\hom(X,-)$ preserves $\kappa$-directed colimits for some $\kappa$ and since $(- \bullet Q)$ is left adjoint, $G$ preserves $\kappa$-directed colimits too. Therefore the colimits of $\kappa$-indexed chains are preserved and due to Theorem \ref{freevar} the corresponding variety has free objects. The rest is obvious.
\end{proo}


\begin{thebibliography}{10}
\bibitem{A74} J. Ad\'{a}mek, {\em Free algebras and automata realizations in the language of categories}, Commentationes Mathematicae Universitatis Carolinae 015 (1974), no. 4, 589-602.
\bibitem{AHS} J. Ad\'{a}mek, H. Herrlich, G. Strecker: {\em Abstract and Concrete categories}, Free Software Foundation (1990/2006)
\bibitem{AP} J. Ad\'{a}mek, H. Porst: {\em From Varieties of Algebras to Coverieties of Coalgebras}, Mathematical Structures in Computer Science (2001)
\bibitem{AR}  J. Ad\'{a}mek, J. Rosick\'{y}: {\em Locally Presentable and Accessible Categories}, Cambridge University Press, (1994)
\bibitem{AT}  J. Ad\'{a}mek, V. Trnkov\'{a}: {\em Birkhoff's Variety Theorem With And Without Free Algebras}, Theory and Applications of Categories, Vol. 14, No. 18, pp 424-450 (2005)
\bibitem{Barr} M. Barr: {\em Coequaliezers and Free Triples}, Mathematische Zeitschrift 116, 307-322, Springer-Verlag (1970)
\bibitem{Kelly} G. M. Kelly: {\em A unified treatment of transfinite constructions for free algebras, free monoids, colimits, associated sheaves, and so on}, Bulletin of the Australian Mathematical Society, 22 (1980)
\bibitem{mac} S. Mac Lane: {\em Categories for the Working Mathematician}, Springer-Verlag (1971)
\bibitem{reit78} J. Reitermann: {\em One More Categorical Model of Universal Algebra}, Mathematische Zeitschrift, Springer-Verlag (1978)
\bibitem{reit86} J. Reitermann: {\em On locally small based algebraic theories}, Commenatationes Mathematicae Unversitatis Carolinae (1986)
\end{thebibliography}
\end{document}